\documentclass[journal]{IEEEtran}

\usepackage{cite}
\usepackage{amssymb,amsmath, amsfonts, amsthm}
\usepackage[caption=false,font=footnotesize]{subfig}
\usepackage{tikz,pgfplots}
\usepackage{pgfplots}
\usetikzlibrary{shapes.geometric, decorations.pathmorphing, positioning, arrows, arrows.meta, patterns, shapes, calc, fit, backgrounds}
\graphicspath{{./figures/AWPL_TDPMCHWT/}}
\DeclareGraphicsExtensions{.pdf,.jpeg,.png,.eps}

\usepackage{multirow}

\newcommand{\R}{\mathbb R}

\newcommand{\NN}{\mathcal N}

\newcommand{\TT}{\mathcal T}
\newcommand{\OO}{\mathcal O}

\newcommand{\KK}{\mathcal K}

\DeclareMathOperator{\xxs}{\mathbf x}

\DeclareMathOperator{\KKs}{\mathbf K}

\DeclareMathOperator{\QQs}{\mathbf Q}
\DeclareMathOperator{\IIs}{\mathbf I}

\DeclareMathOperator{\diag}{diag}

\DeclareMathOperator{\TTs}{\mathbf T}
\DeclareMathOperator{\ZZs}{\mathbf Z}

\newcommand{\paren}[1]{\left({#1}\right)}

\newcommand{\abs}[1]{\left\vert{#1}\right\vert}

\newcommand{\transpose}{\mathsf{T}}

\newcommand{\bm}[1]{\boldsymbol{#1}}

\newcommand{\curlt}{\textbf{\textup{curl}}}
\newcommand{\divt}{\textup{div}}
\newcommand{\gradt}{\textbf{grad}}

\newcommand{\veps}{\varepsilon}

\newcommand{\pa}{\partial}

\newcommand{\Gm}{\Gamma}

\newcommand{\sm}{\sigma}

\newcommand{\xb}{\bm{x}}
\newcommand{\yb}{\bm{y}}

\newcommand{\hb}{\bm{h}}

\newcommand{\eb}{\bm{e}}

\newcommand{\ub}{\bm{u}}

\newcommand{\fb}{\bm{f}}
\newcommand{\gb}{\bm{g}}

\newcommand{\jb}{\bm{j}}

\newcommand{\mb}{\bm{m}}

\newcommand{\nv}{\textbf{n}}
\newcommand{\cv}{\textbf{c}}
\newcommand{\uv}{\textbf{u}}
\newcommand{\jv}{\textbf{j}}
\newcommand{\mv}{\textbf{m}}
\newcommand{\rv}{\textbf{r}}
\newcommand{\ev}{\textbf{e}}
\newcommand{\hv}{\textbf{h}}

\newcommand{\zrb}{\bm{0}}

\newcommand{\ds}{\, \mathrm{d}s}

\newcommand{\di}{\, \mathrm{d}}

\newcommand{\q}{\quad}

\newcommand{\qqq}{\qquad\quad}
\newcommand{\qqqq}{\qquad\qquad}
\newcommand{\qqqqq}{\qquad\qquad\quad}

\begin{document}

\title{On the Late-Time Instability of MOT solution to the Time-Domain PMCHWT Equation}

\author{Van~Chien~Le,~\IEEEmembership{Member,~IEEE,}
        Viviana~Giunzioni,~\IEEEmembership{Graduate Student Member,~IEEE,}
        Pierrick~Cordel,~\IEEEmembership{Graduate Student Member,~IEEE,}
        Francesco~P.~Andriulli,~\IEEEmembership{Fellow,~IEEE,}
        and~Kristof~Cools,~\IEEEmembership{Member,~IEEE}
\thanks{Manuscript received April 19, 2005; revised August 26, 2015. This work was supported by the European Research Council (ERC) under the European Union’s Horizon 2020 research and innovation programme (Grant agreement No. 101001847).}%
\thanks{V.~C. Le and K. Cools are with IDLab, Department of Information Technology at Ghent University - imec, 9000 Ghent, Belgium (e-mail: vanchien.le@ugent.be, kristof.cools@ugent.be).}
\thanks{V. Giunzioni, P. Cordel and F.~P. Andriulli are with the Department of Electronics and Telecommunications, Politecnico di Torino, 10129 Turin, Italy.}}

\markboth{IEEE Antennas and Wireless Propagation Letters}{LE \MakeLowercase{\textit{et al.}}: On the late-time instability of the TD-PMCHWT equation}

\maketitle

\begin{abstract}
    This paper investigates the late-time instability of marching-on-in-time solution to the time-domain PMCHWT equation. The stability analysis identifies the static solenoidal nullspace of the time-domain electric field integral operator as the primary cause of instability. Furthermore, it reveals that the instability mechanisms of the time-domain PMCHWT equation are fundamentally different from those of the time-domain electric field integral equation. In particular, the PMCHWT's instability is much more sensitive to numerical quadrature errors, and its spectral characteristics are strongly influenced by the topology and smoothness of the scatterer surface.
\end{abstract}



\section{Introduction}

\IEEEPARstart{T}{he} time-domain Poggio-Miller-Chang-Harrington-Wu-Tsai (TD-PMCHWT) equation is among the most common formulations to model transient scattering by (piecewise) homogeneous dielectric bodies \cite{PM1973}. Besides the ill-conditioned nature, the TD-PMCHWT equation is particularly susceptible to late-time instability. The instability mechanisms of the TD-PMCHWT are not nearly as well understood as is the case for the time-domain electric field integral equation (TD-EFIE). Unlike the TD-EFIE operator (TD-EFIO) that supports the static solenoidal nullspace, the continuous TD-PMCHWT operator does not possess a nullspace. Nevertheless, in marching-on-in-time (MOT) schemes, instability emerges, superficially resembling direct-current (DC) instability of the TD-EFIE, but much more sensitive to numerical quadrature errors \cite{MB1982,BCA2015d,Beghein2015,LAC2024}.

Several techniques have been proposed to mitigate this issue. Most of them concern improving the accuracy of integral evaluations, including averaging/filtering techniques \cite{SA1993,Rynne1994,DD1997}, the use of appropriate functions \cite{PBW1998,WPC+2004,PW2006,YT2017} and accurate integration schemes \cite{YE2006,SLY+2009,SXC+2011,VVV+2013,VAB+2013,TXX2014}. While these methods provide partial improvements, they are effective only in specific cases or to a limited extent (at most a constant solution at late times). A comprehensive understanding of the mechanisms of instability is crucial for developing robust stabilization strategies. 

In this letter, we investigate the late-time instability of MOT solution to the TD-PMCHWT equation through a stability analysis based on the spectral properties of a companion matrix. The analysis shows that the number of unstable modes is linked to the TD-EFIO. However, it also reveals that the behavior of the TD-PMCHWT instability is fundamentally different from that of the TD-EFIE. In particular, we quantitatively characterize the late-time instability in terms of numerical quadrature and round-off errors, highlighting the sensitivity of TD-PMCHWT to the accuracy of interaction matrix computations. Furthermore, we demonstrate that the spectral characteristics of the instability are significantly influenced by geometric properties of the scatterer surface.

\section{Marching-on-in-Time Scheme}
\label{sec:formulation}

Let $\Gm$ be the surface of a body in $\R^3$ filled by a homogeneous dielectric material with permittivity $\epsilon^\prime$ and permeability $\mu^\prime$. This body is immersed in a homogeneous background medium $(\epsilon, \mu)$. An incident transient wave $(\eb^{in}, \hb^{in})$ induces on $\Gm$ the electric and magnetic current densities $\jb(\xb, t)$ and $\mb(\xb, t)$, which satisfy the TD-PMCHWT equation
\begin{equation}
\label{eq:PMCHWT}
    \begin{pmatrix}
        \eta \TT + \eta^\prime \TT^\prime  & -\KK - \KK^\prime \\
        \KK + \KK^\prime & \frac{1}{\eta} \TT + \frac{1}{\eta^\prime} \TT^\prime
    \end{pmatrix}
    \begin{pmatrix}
        \jb \\
        \mb
    \end{pmatrix} 
    =
    \begin{pmatrix}
        \eb^{in} \times \nv \\
        \hb^{in} \times \nv
    \end{pmatrix}.
\end{equation}
 The TD-EFIO $\TT$ and the time-domain magnetic field integral operator (TD-MFIO) $\KK$ associated with $(\epsilon, \mu)$ are given by 
\begin{align*}
    (\TT \jb)(\xb, t) & = (\TT^s \jb)(\xb, t) + (\TT^h \jb)(\xb, t), \\  
    (\TT^s \jb)(\xb, t) & = - \dfrac{1}{c} \, \nv \times \int_{\Gm} \dfrac{\pa_t \jb(\yb, \tau)}{4\pi R} \ds_{\yb}, \\
    (\TT^h \jb)(\xb, t) & = c \, \nv \times \gradt_{\xxs} \int_{\Gm} \int_{-\infty}^{\tau}\dfrac{\divt_\Gm \, \jb(\yb, \xi)}{4\pi R} \di \xi \ds_{\yb}, \\
    (\KK\jb)(\xb, t) & = \nv \times \curlt_{\xxs} \int_{\Gm} \dfrac{\jb(\yb, \tau)}{4\pi R} \ds_{\yb}.
\end{align*}
Here, $\eta = \sqrt{\mu/\epsilon}, c = 1/\sqrt{\mu\epsilon}$, $R = \abs{\xb - \yb}, \tau = t - R/c$, and $\nv$ is the outward normal of $\Gm$. The operators and quantities $(\TT^\prime, \KK^\prime, \eta^\prime, c^\prime)$ associated with $(\epsilon^\prime, \mu^\prime)$ are defined analogously.

For comparison, we recall the TD-EFIE, which models the surface electric current density $\jb(\xb, t)$ induced by an incident wave at the perfect electric conductor surface $\Gm$ as
\[
    \eta \TT \jb = \eb^{in} \times \nv.
\]

It is well-known that the TD-EFIOs $\TT$ and $\TT^\prime$ have a common nullspace comprising all static solenoidal currents \cite{ACO+2009}. In contrast, the TD-PMCHWT operator does not possess any nullspace (one can easily verify for the case $(\epsilon^\prime, \mu^\prime) = (\epsilon, \mu)$ using the time-domain Calder\'{o}n identities in \cite{CAO+2009b}). Furthermore, the TD-EFIE's solution is affected by instability arising from resonant frequencies of the interior body \cite{LCA+2024}, whereas the TD-PMCHWT's solution is free from that \cite{BHV+2003}.

The boundary $\Gm$ is partitioned into a mesh of $N_f$ triangles with $N_e$ edges, equipped with the Rao-Wilton-Glisson (RWG) basis functions $\fb_n(\xb)$, with $n = 1, 2, \ldots, N_e$ \cite{RWG1982}. Along the time axis, $N_t$ time intervals of length $\Delta t$ are equipped with the hat basis functions $h_i(t)$, with $i = 1, 2, \ldots, N_t$.

The current densities $\jb$ and $\mb$ are approximated by their expansions in the trial space via the coefficient vectors $\jv$ and $\mv$, respectively. The TD-PMCHWT \eqref{eq:PMCHWT} is tested with the spatial rotated RWG functions and temporal Dirac delta distributions, resulting in a lower-triangular block matrix system that can be solved using the MOT scheme
\begin{equation}
\label{eq:MOT}
    \uv_i = \ZZs_0^{-1} \paren{\rv_i - \sum_{k=1}^{i-1} \ZZs_k \uv_{i-k}}, \q \,\,\, i = 1, 2, \ldots, N_t,
\end{equation}
where $\uv_i = \paren{\jv_i, \mv_i}^{\transpose}, \rv_i = \paren{\ev_i, \hv_i}^{\transpose}$, and
\[
    \ZZs_k  = 
    \begin{pmatrix}
        \eta \TTs_k + \eta^\prime \TTs_k^\prime  & -\KKs_k - \KKs_k^\prime \\
        \KKs_k + \KKs_k^\prime & \frac{1}{\eta} \TTs_k + \frac{1}{\eta^\prime} \TTs_k^\prime
    \end{pmatrix}.
\]
The MOT system for the TD-EFIE is derived analogously. In most cases, it is not necessary to evaluate all matrices $\ZZs_k$ due to the finite support of the temporal basis functions $h_i(t)$ \cite{VDV2023}. More specifically, for all $k > k_0 := \left\lceil \tfrac{D}{c_{\text{min}} \Delta t} \right\rceil$, with $D$ the diameter of the scatterer and $c_{\text{min}} = \min(c, c^\prime)$, the contributions of $\TT^s, \KK$ and their inner counterparts to $\ZZs_k$ vanish, while the contributions of $\TT^h$ and its inner counterpart are given by
\[
    \left[\TTs_\infty^h\right]_{mn} = -\Delta t \int_\Gm \int_\Gm \dfrac{\divt_{\Gm} \fb_m(\xb) \, \divt_{\Gm} \fb_n(\yb)}{4\pi R} \ds_{\yb} \ds_{\xb}.
\]
In particular, for all $k > k_0$,
\[
    \ZZs_k = \ZZs_{\infty} := \diag\paren{\dfrac{\epsilon + \epsilon^\prime}{\epsilon \epsilon^\prime} \TTs_\infty^h, \dfrac{\mu + \mu^\prime}{\mu \mu^\prime} \TTs_\infty^h}.
\]

\begin{figure}[t!]
\begin{tikzpicture}
    \node[anchor=north west] at (0,0) {\includegraphics[width=\linewidth]{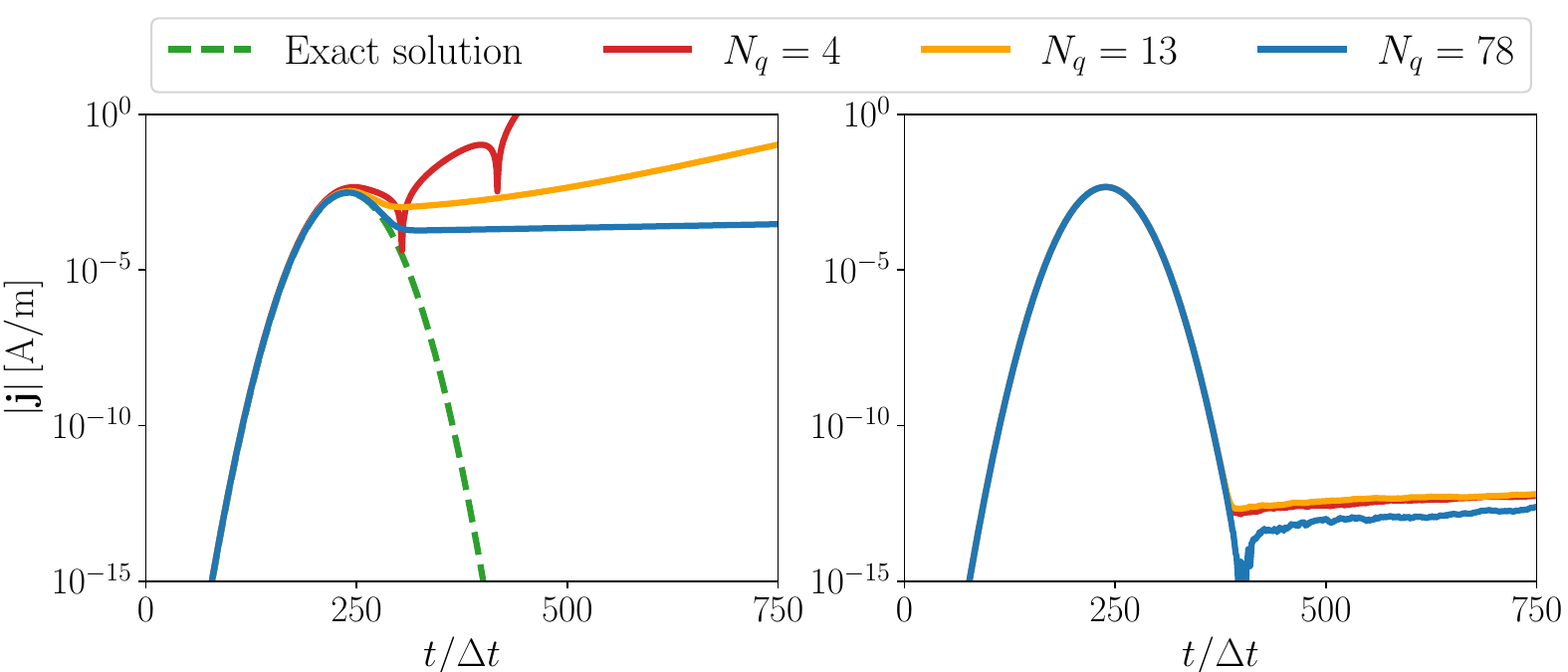}};
    \node[anchor=north west] at (7,-1) {\includegraphics[trim={36cm 14cm 35cm 16cm}, clip, width=0.17\linewidth]{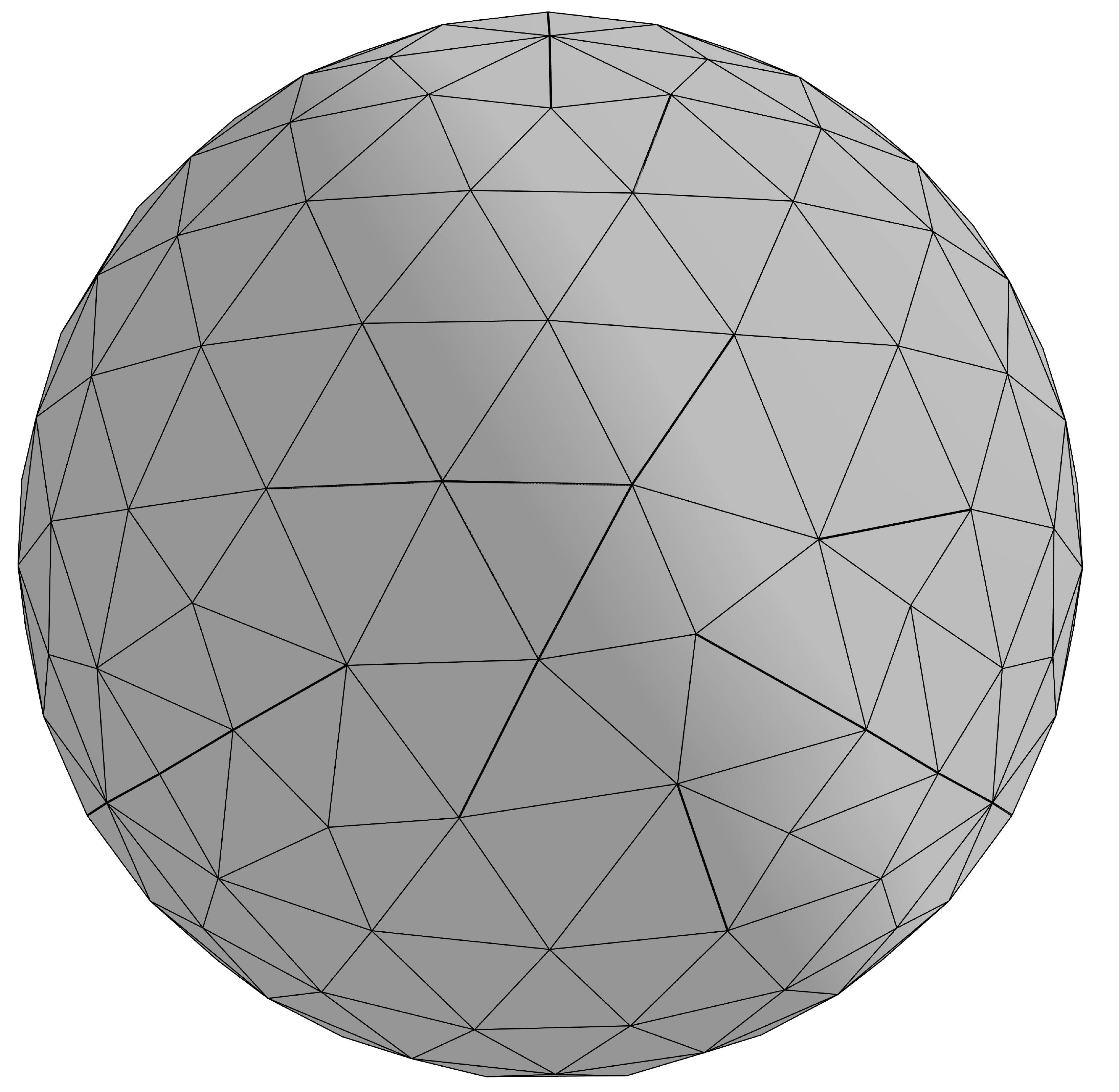}};
    \label{fig:sphere_current}
\end{tikzpicture}
\caption{MOT solutions to the TD-PMCHWT equation (\textit{left}) and the TD-EFIE (\textit{right}) on a sphere. Both solutions suffer from late-time instability. The TD-PMCHWT's instability is significantly affected by the number of quadrature points $N_q$ used to evaluate integrals, whereas that of the TD-EFIE is unaffected.}
\label{fig:solution}
\end{figure}

MOT solutions to the TD-EFIE and TD-PMCHWT equation both suffer from late-time instability, which manifests itself in non-decaying or even exponentially growing errors over time. Fig.~\ref{fig:solution} demonstrates the late-time instability of MOT solutions and the impact of numerical quadrature errors on instability. Here and throughout the paper, a semi-analytic quadrature method is used to discretize the time-domain boundary integral operators, which evaluates the inner integral based on the Wilton technique \cite{WRG+1984,SLY+2009}, while the outer (test) integral is done numerically using $N_q$ quadrature points. The numerical result shows that the level of DC instability (i.e., the exponential growth rate of solution) of the TD-EFIE is typically very low and unaffected by numerical quadrature errors. In contrast, the TD-PMCHWT's instability is very sensitive to the accuracy of evaluating interaction matrices. This difference will be explained in the next section.

\section{Companion Matrix}
\label{sec:companion}

In order to analyze the late-time instability of MOT schemes, we employ the commonly used companion-matrix stability analysis (CMSA) technique \cite{DWB1998}. The reader is also referred to the positive-definite stability analysis technique for a lower computational complexity, especially when MOT systems involve a large number of spatial unknowns \cite{VVD2024}. The CMSA qualitatively determines stability of a MOT scheme by examining whether its solution remains bounded over an unbounded interval, by means of classical stability analysis of numerical methods for ODEs \cite{Butcher2008}. It is done based on the spectral behavior of a companion matrix relating successive sets of solution vectors when the incident field has passed. 

To define the companion matrix for the MOT system \eqref{eq:MOT}, we set $\rv_i = \zrb$ and introduce the composite current vector 
\[
    \textstyle
    \cv_i := \paren{\uv_i, \uv_{i-1}, \ldots, \uv_{i-k_0+1}, \sum_{k=1}^{i-k_0} \uv_k}^\transpose.
\]
Two arbitrary consecutive vectors $\cv_i$ satisfy the recurrence relation $\cv_{i} = \QQs \cv_{i-1}$, where the companion matrix $\QQs$ is given by
\[
    \QQs = 
    \begin{pmatrix}
        \QQs_1 & \QQs_2 & \cdots & \cdots & \QQs_{k_0} & \QQs_{\infty} \\[0.15cm]
        \IIs & \zrb & \cdots & \cdots & \zrb & \zrb \\
        \zrb & \IIs & \ddots & \ddots & \vdots & \vdots \\
        \vdots & \vdots & \ddots & \ddots & \vdots & \vdots \\
        \zrb & \zrb & \cdots & \IIs & \zrb & \zrb \\[0.15cm]
        \zrb & \zrb & \cdots & \zrb & \IIs & \IIs
    \end{pmatrix},
\]
with $\QQs_k := -\ZZs_0^{-1} \ZZs_k, k = 1, 2, \ldots, k_0, \infty$, and $\IIs$ the identity matrix of size $2N_e \times 2N_e$ \cite{VDZ+2022}.

Next, we investigate the spectrum of the companion matrix $\QQs$ for the TD-EFIE and TD-PMCHWT. Let $\NN_\Gm$ be the space spanned by RWG coefficient vectors of solenoidal functions, which has dimension $\dim \NN_\Gm = N_l := N_e - N_f + 1$ \cite{Andriulli2012}. Obviously, $\NN_\Gm$ is the nullspace of the matrix $\TTs_\infty^h$. As a consequence, the companion matrix associated with the TD-EFIE has $N_l$ ``trivial'' eigenvalues $\lambda = 1 + 0j$, which are corresponding to the ordinary eigenvectors
\[
    \cv_0 = \paren{\zrb, \zrb, \ldots, \zrb, \gb_h}^\transpose, \qqqqq \gb_h \in \NN_\Gm.
\]
In addition, the matrix $\QQs$ of the TD-EFIE has other $N_l$ ``non-trivial'' eigenvalues $\lambda = 1+0j$ associated with the static solenoidal nullspace of the TD-EFIO \cite{ACO+2009}. Precisely, they are corresponding to the generalized eigenvectors of rank 2 
\[
    \tilde{\cv}_0 = \paren{\gb_h, \gb_h, \ldots, \gb_h, \zrb}^\transpose, \qqqq \gb_h \in \NN_\Gm.
\]

Similarly to the TD-EFIE, the companion matrix $\QQs$ of the TD-PMCHWT has $2 N_l$ trivial eigenvalues $\lambda = 1 + 0j$, which are associated with the nullspace of the tail $\ZZs_\infty$. In contrast, the non-trivial eigenvalues of the TD-EFIE do not belong to the spectrum of the TD-PMCHWT. Nevertheless, in numerical experiments, the spectrum of the TD-PMCHWT reveals additional eigenvalues clustered around $1 + 0j$, apart from the $2N_l$ trivial ones. These nearby eigenvalues, while not strictly equal to $1 + 0j$, are referred to as the non-trivial eigenvalues of the TD-PMCHWT. 

Upon discretization, numerical errors arise, primarily comprising dominant quadrature errors $\sm$ and round-off errors $\veps_{\text{mach}}$ ($\approx 10^{-15}$ in double precision). These errors cause shifts in the eigenvalues of the companion matrix. In particular, eigenvalues originally located at or near $1 + 0j$ that are shifted outside the unit circle lead to the exponential growth of MOT solution $\uv_i \sim (1+r)^i$ at late times, with growth rate $r$ of the order of the maximum eigenvalue shift $\delta$.

Since the eigenvectors $\cv_0$ and $\tilde{\cv}_0$ can be discretized with machine precision $\veps_{\text{mach}}$ and the $2N_l$ eigenvalues at $1 + 0j$ of the TD-EFIE are coupled into $N_l$ Jordan blocks of size 2, these eigenvalues are symmetrically shifted by $\delta = \OO(\sqrt{\veps_{\text{mach}}})$\cite{VDZ+2022}. In contrast, the trivial and non-trivial eigenvalues of the TD-PMCHWT are uncoupled. The trivial eigenvalues are affected by small shifts $\OO(\veps_{\text{mach}})$, while the non-trivial ones are perturbed by the quadrature errors $\OO(\sm)$ \cite{Golub2013}. Consequently, the MOT solution to the TD-EFIE exhibits an exponential growth at a rate $\OO(\sqrt{\veps_{\text{mach}}})$, whereas the TD-PMCHWT solution grows at a rate $\OO(\sm)$.

\begin{figure}[!t]
    \centering
    \includegraphics[trim={0.2cm 0cm 0.5cm 0.5cm}, clip, width=\linewidth]{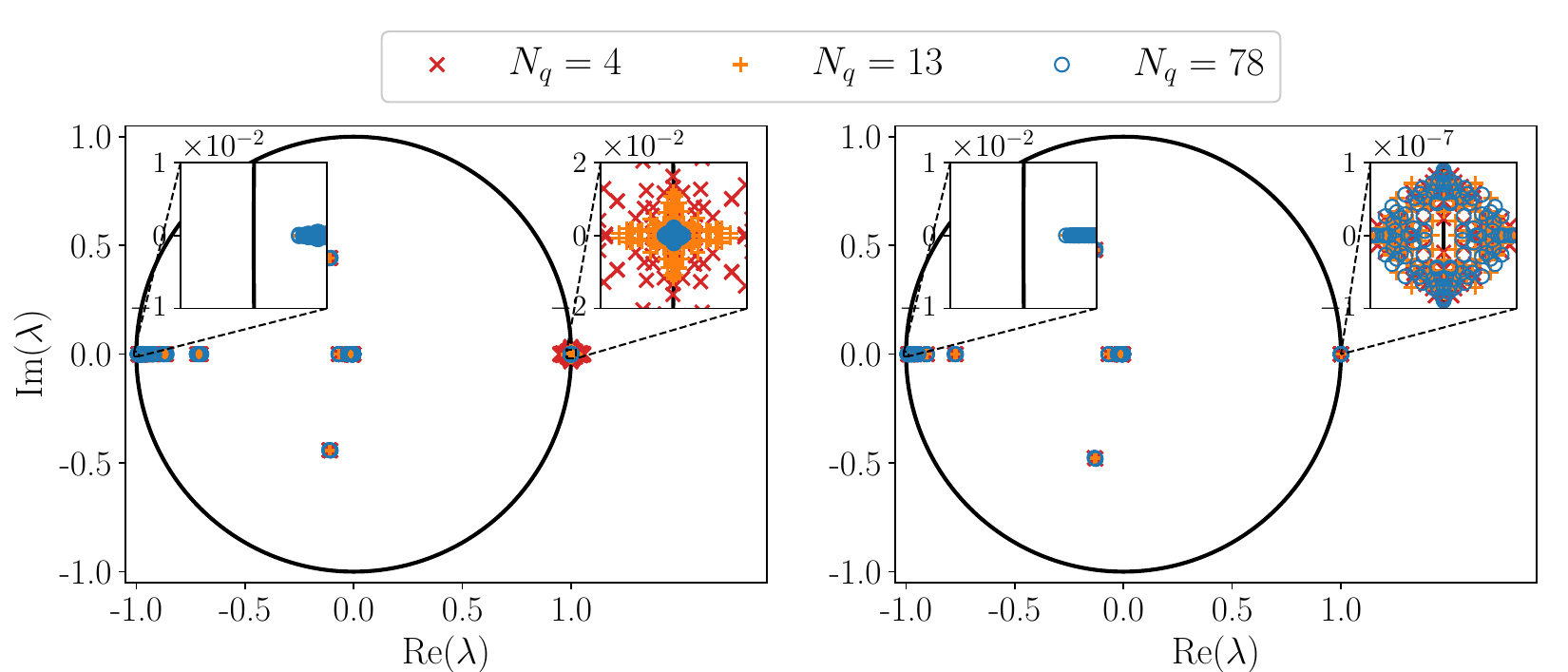}
    \caption{Spectrum of the companion matrices for the TD-PMCHWT (\textit{left}) and the TD-EFIE (\textit{right}) on a sphere. The shift $\delta$ of eigenvalues near $1 + 0j$ of the TD-PMCHWT decreases when increasing the number of quadrature points $N_q$, whereas that of the TD-EFIE is independent of quadrature rules (see Table~\ref{tab:error_shift} for further details). The eigenvalues near $-1 + 0j$ are not problematic as they reside strictly inside the unit circle \cite{VVV+2013}.}
\label{fig:sphere_eigen}
\end{figure}

This fundamental difference in late-time instability of the two formulations is corroborated by Fig.~\ref{fig:sphere_eigen} and Table~\ref{tab:error_shift}, and is further illustrated in the solution behavior shown in Fig.~\ref{fig:solution}. These findings provide a rigorous explanation for observations in \cite{SLY+2009,Beghein2015,LAC2024}, where the TD-PMCHWT solution was found to be more susceptible to instability than the TD-EFIE.

In numerical experiments, the eigenvalue shift $\delta$ and the quadrature error $\sm$ are estimated by
\[
    \delta = \abs{\lambda_{\max}(\QQs)} - 1, \qqq \sm = \abs{\lambda_{\max}(\QQs - \QQs_{ref})},
\]
where $\QQs_{ref}$ is the reference companion matrix computed with a very accurate quadrature rule ($N_q = 400$).

\setlength{\arrayrulewidth}{0.2mm}
\setlength{\tabcolsep}{18pt}
\renewcommand{\arraystretch}{1.28}

\begin{table}[!t]
\centering
\caption{Quadrature errors $\sm$ vs. eigenvalue shifts $\delta$}
\label{tab:error_shift}
\setlength{\tabcolsep}{5pt}
\begin{tabular}{|c|c|c|c|c|}
\hline
\multirow{2}{*}{} & \multicolumn{2}{|c|}{\bf TD-PMCHWT} & \multicolumn{2}{|c|}{\bf TD-EFIE} \\
\cline{2-5}
 & \textbf{Error} $\sm$ & \textbf{Shift} $\delta$ & \textbf{Error} $\sm$ & \textbf{Shift} $\delta$ \\
\hline
$N_q = 4$ & $5.26 \cdot 10^{-2}$ & $5.74 \cdot 10^{-2}$ & $3.56 \cdot 10^{-3}$ & $1.06 \cdot 10^{-7}$ \\ \hline
$N_q = 13$ & $1.66 \cdot 10^{-2}$ & $1.55 \cdot 10^{-2}$ & $1.28 \cdot 10^{-4}$ & $1.24 \cdot 10^{-7}$ \\
\hline
$N_q = 78$ & $3.11 \cdot 10^{-3}$ & $2.52 \cdot 10^{-3}$ & $1.56 \cdot 10^{-6}$ & $1.14 \cdot 10^{-7}$ \\
\hline
\end{tabular}
\end{table}

\section{Stability Analysis}
\label{sec:analysis}

In this section, the CMSA technique is applied to study the origin of the non-trivial eigenvalues of the TD-PMCHWT, as well as the effects of different factors on late-time instability.

\subsection{Origin of non-trivial eigenvalues and effects of geometry}

\subsubsection{Simply-connected domains}

\begin{figure}[!t]
    \centering
    \includegraphics[trim={2cm 3cm 3cm 2cm}, clip, width=0.15\linewidth]{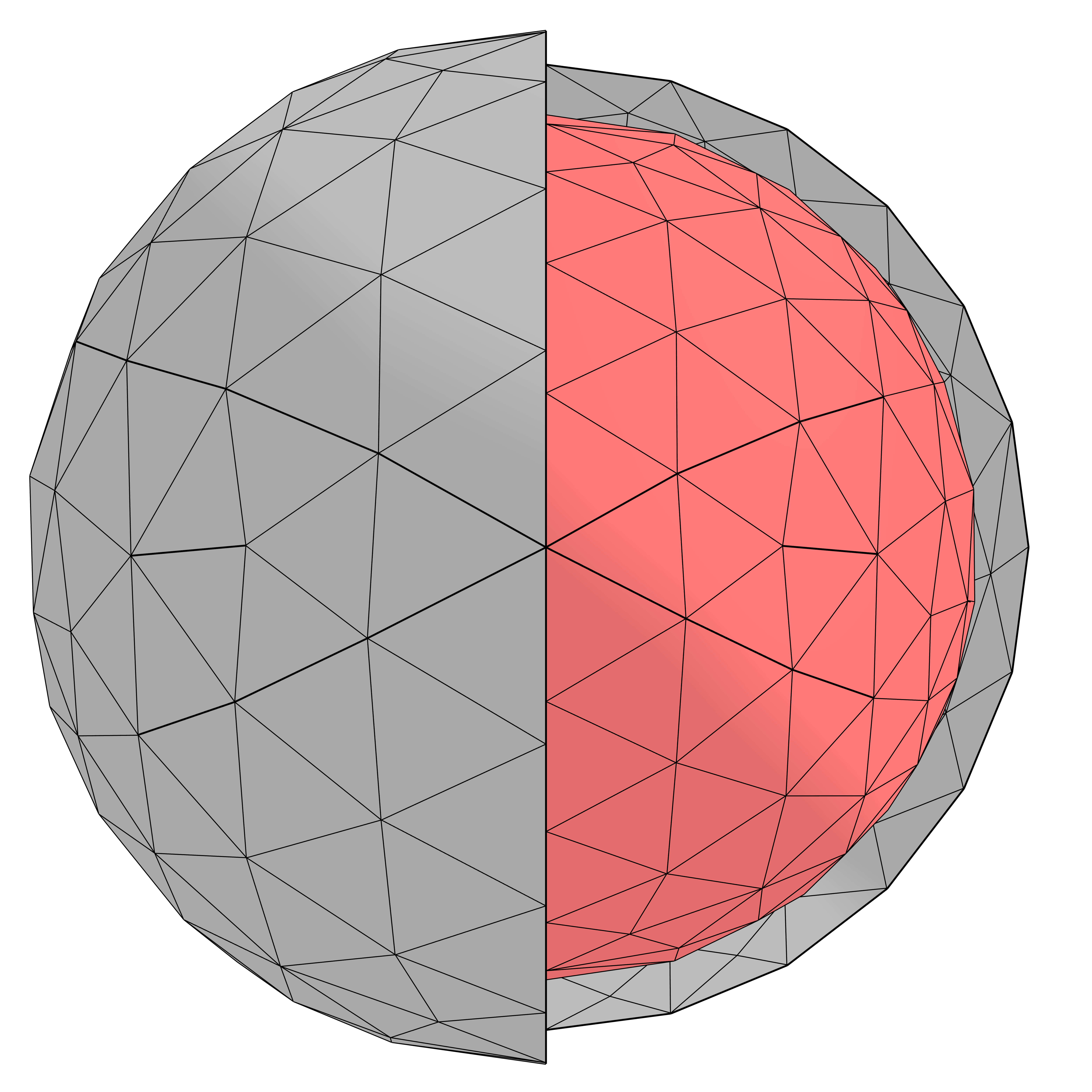}
    \hfill
    \includegraphics[trim={29cm 10cm 29cm 25cm}, clip, width=0.21\linewidth]{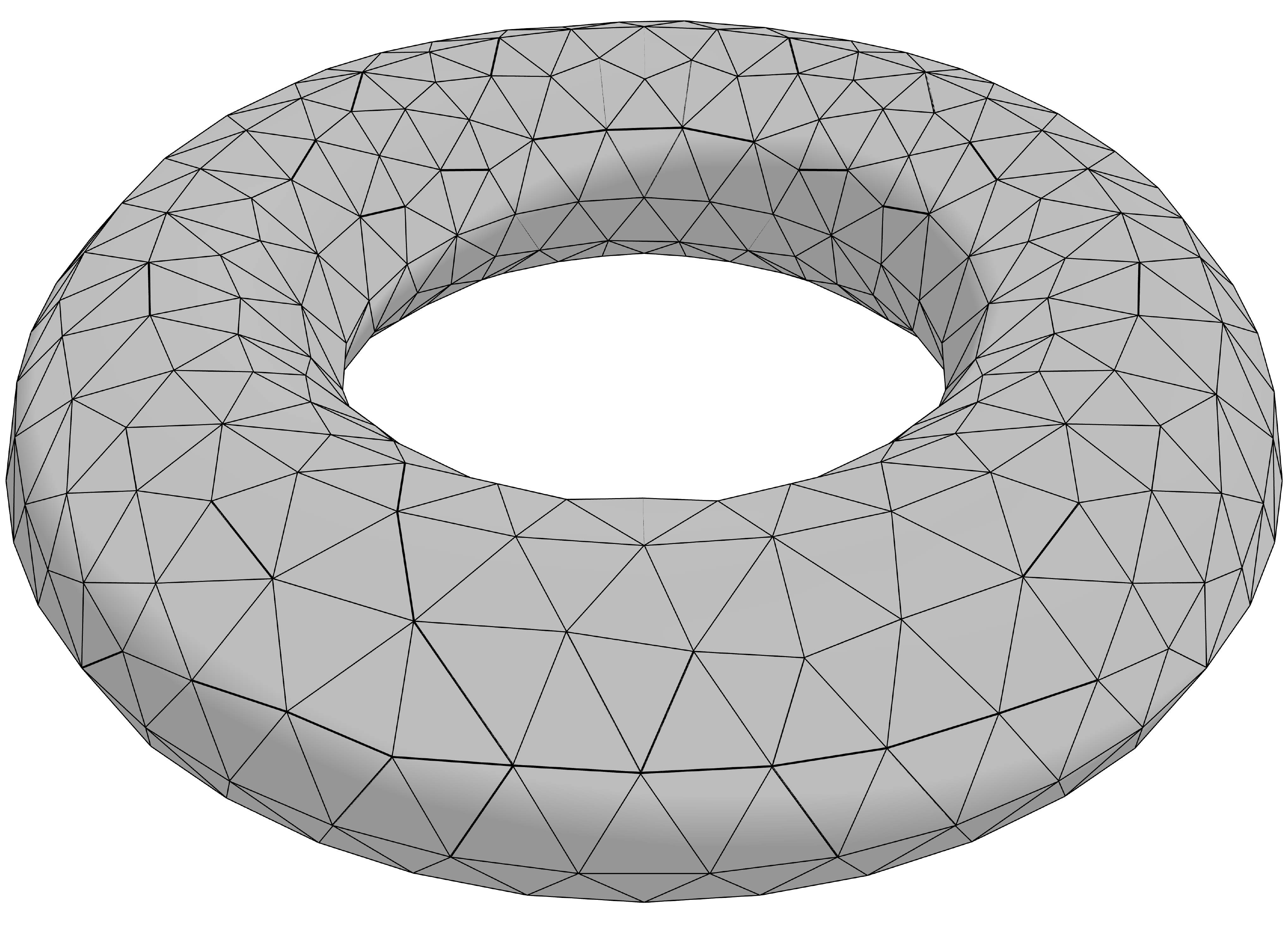}
    \includegraphics[trim={12cm 13cm 23cm 22cm}, clip, width=0.36\linewidth]{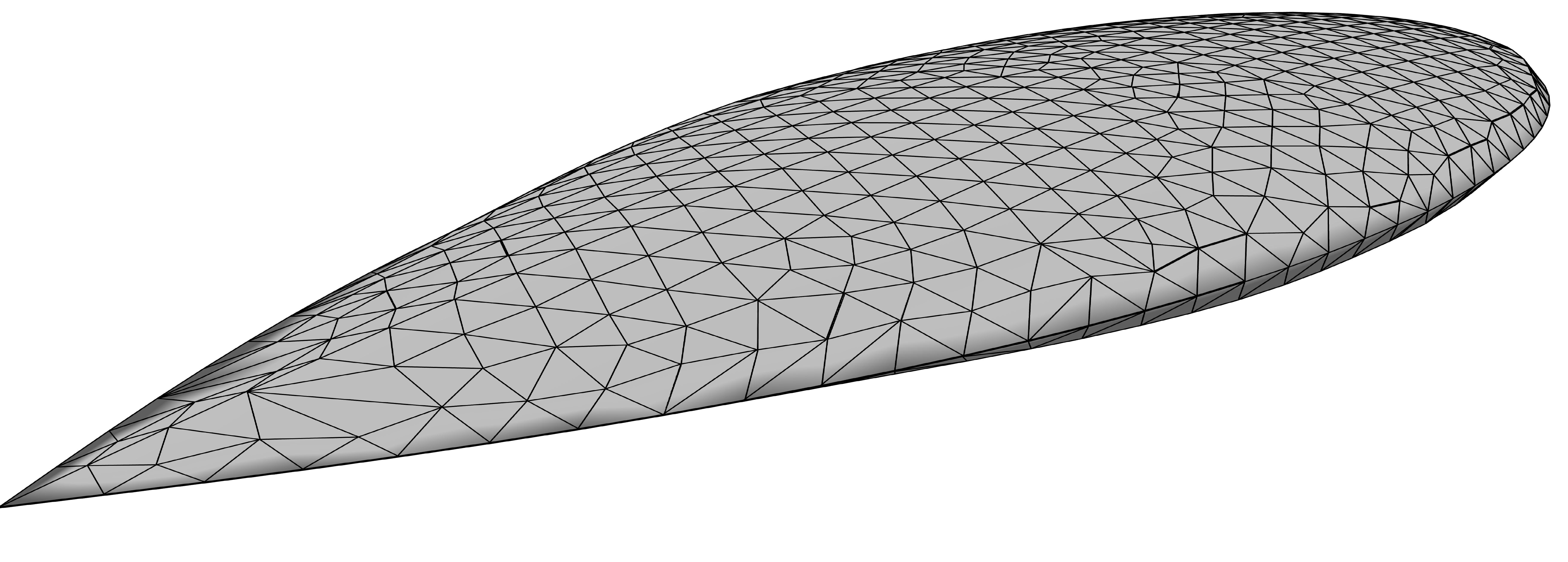}
    \hfill
    \includegraphics[trim={34cm 26cm 34cm 18cm}, clip, width=0.23\linewidth]{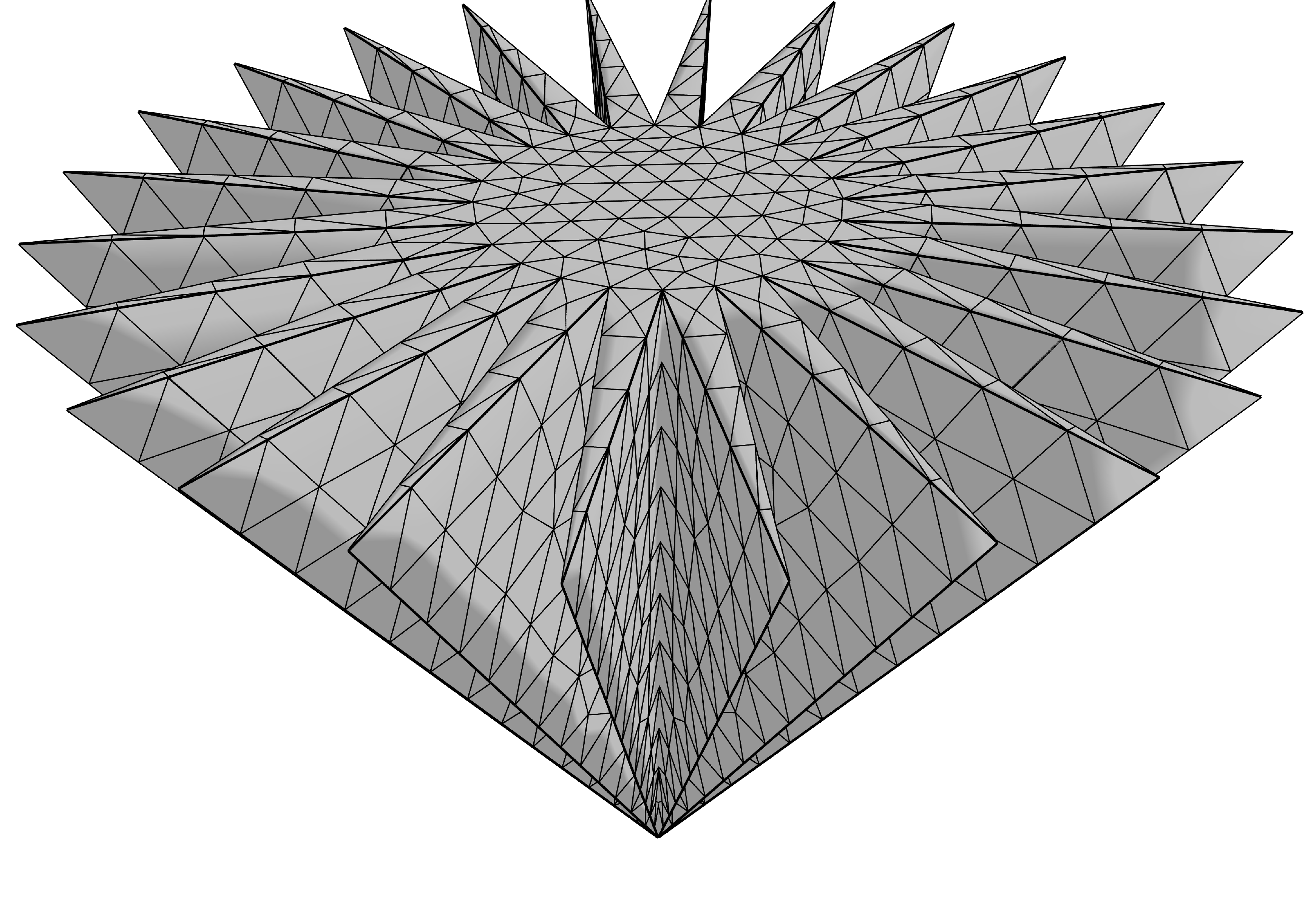}
    \caption{Geometries exhibiting different types of challenges. From left to right: a two-layered sphere of radii $1\mathrm{m}$ and $0.85\mathrm{m}$; a torus of radii $0.75\mathrm{m}$ and $0.25\mathrm{m}$; a NASA almond of size $9.936 \mathrm{m} \times 1.92 \mathrm{m} \times 0.64 \mathrm{m}$ \cite{WWS+1993}; a pyramid of height $0.5 \mathrm{m}$ and 24-pointed star base, whose vertices lie on two concentric circles of radius $1\mathrm{m}$ and $0.3\mathrm{m}$ \cite{LC2024b}.}
    \label{fig:geometries}
\end{figure}

\begin{figure*}[!t]
\centering
\subfloat{
    \includegraphics[trim={0.3cm 0cm 0.5cm 0.9cm}, clip, width=0.48\linewidth]{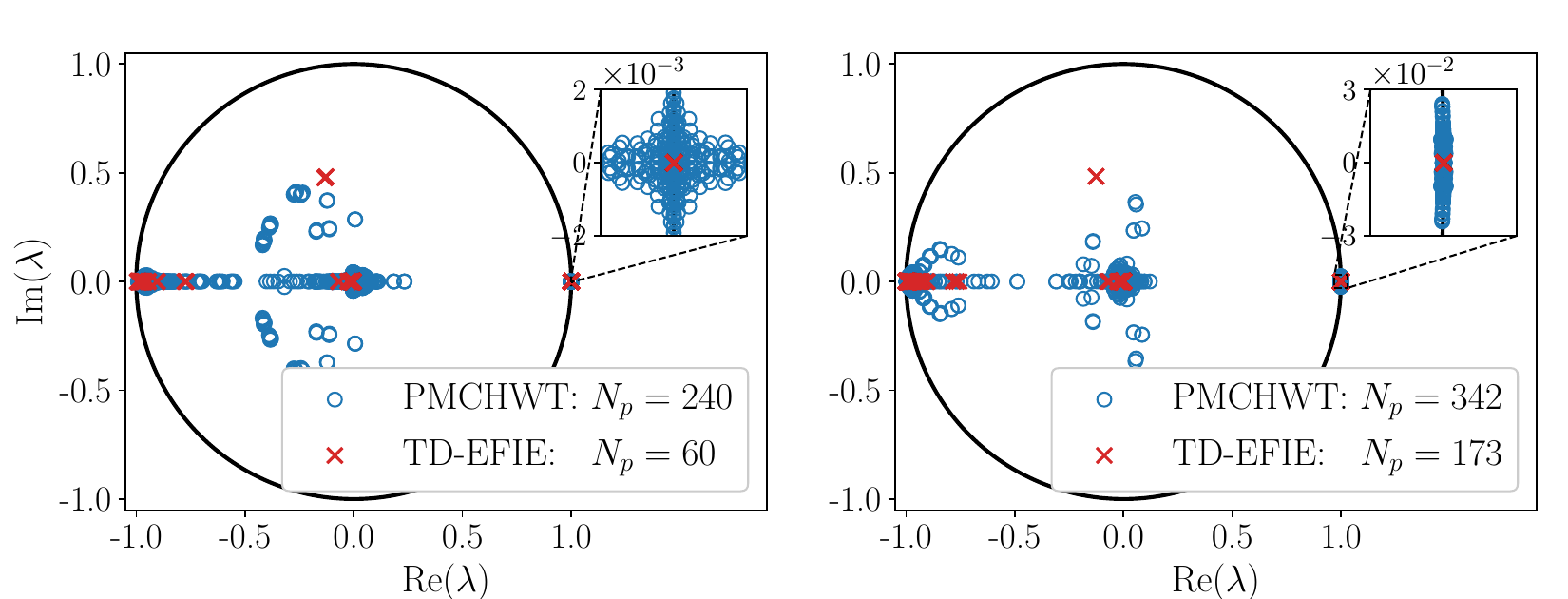}
    }
\hfil
\subfloat{
    \includegraphics[trim={0.3cm 0cm 0.5cm 0.9cm}, clip, width=0.48\linewidth]{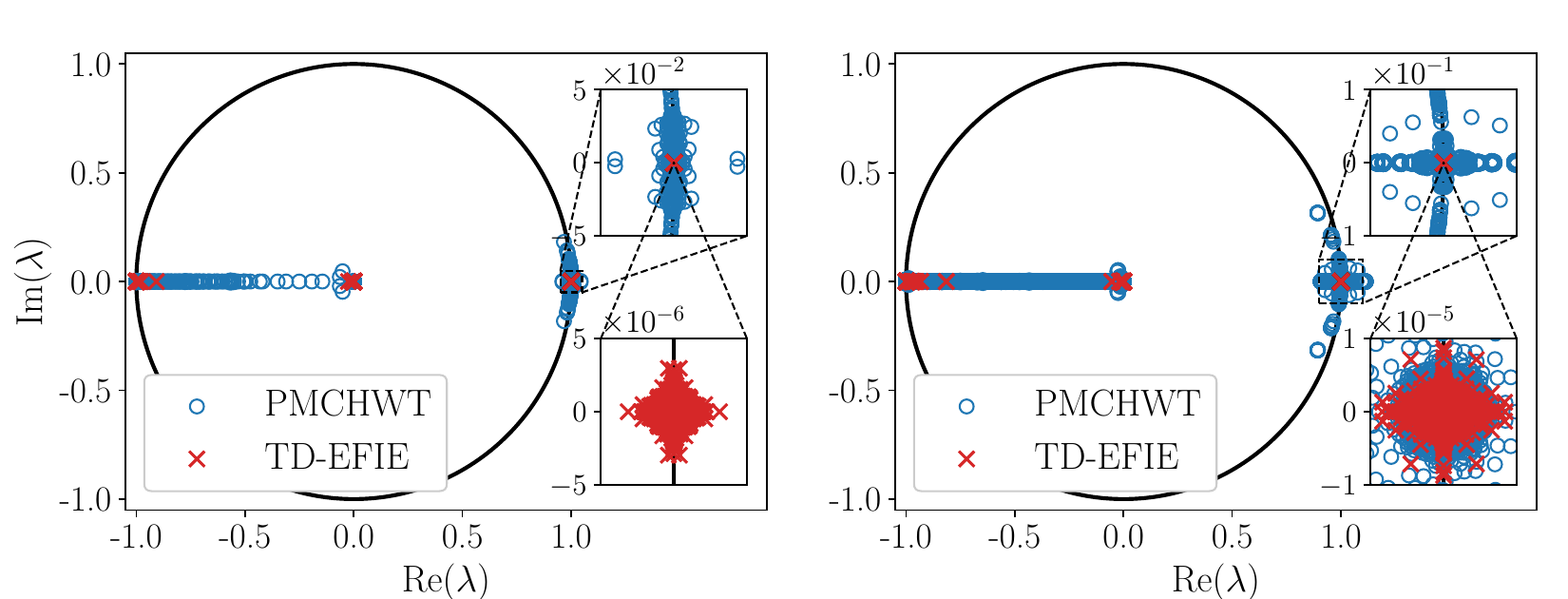}
    }
\caption{Spectrum of the companion matrices for the TD-PMCHWT and the TD-EFIE on a two-layered sphere, a torus, a NASA almond, and a star-based pyramid (from left to right). There are $N_p = 2 N_l$ non-trivial eigenvalues of the TD-PMCHWT on the sphere, while only $N_p = 2N_l - 4$ non-trivial eigenvalues cluster at $1+0j$ for the torus. On the non-smooth domains, the eigenvalues around $1 + 0j$ of the TD-EFIE remain symmetrically distributed, while those of the TD-PMCHWT largely spread to the interior of the unit circle, reflecting the significant impact of the TD-MFIOs on the spectral distribution.}
\label{fig:eigen}
\end{figure*}

We first perform a stability analysis of the MOT system for a two-layered sphere (see Fig.~\ref{fig:geometries} and Table~\ref{tab:shift_geometry}). The spectrum of the companion matrices for this simply-connected smooth surface is depicted in Fig.~\ref{fig:eigen} (\textit{left}). The PMCHWT spectrum exhibits $N_p = 2N_l$ non-trivial eigenvalues clustered around $1+0i$, in addition to $2N_l$ trivial ones. This behavior indicates that the discrete TD-PMCHWT operator on the sphere inherits the static solenoidal nullspace of the TD-EFIOs. For smooth geometries, the MFIOs are known to be compact operators \cite{Nedelec2001,BH2003}. As a result, the contributions of $\KK$ and $\KK^\prime$, along with numerical errors, can be viewed as perturbations to $\TT$ and $\TT^\prime$. 

This analysis yields a key conclusion: the non-trivial eigenvalues near $1 + 0j$ of the TD-PMCHWT are primarily linked to the static solenoidal nullspace of the TD-EFIOs. This nullspace is therefore the principal source of instability. On smooth simply-connected surfaces, the TD-MFIOs and quadrature errors mainly act as numerical perturbations.

\renewcommand{\arraystretch}{1.2}

\begin{table}[!t]
\footnotesize
\centering
\caption{Eigenvalue shifts $\delta$ vs. geometries}
\label{tab:shift_geometry}
\setlength{\tabcolsep}{5pt}
    \begin{tabular}{|c|c|c|c|c|}
        \hline
        \multirow{2}{*}{\bf Geometry} & \multirow{2}{*}{\bf Challenges} & \multicolumn{3}{|c|}{\bf Eigenvalue shift $\delta$} \\
        \cline{3-5}
        & & $N_q = 4$ & $N_q = 13$ & $N_q = 78$ \\ \hline
        Layered sphere  & Multi-domain & $0.0595$ & $0.0133$ & $0.0022$ \\ \hline
        Torus  & Multiply-connected & $0.0118$ & $0.0045$ & $0.0009$ \\ \hline
        NASA almond & Highly non-smooth & $0.7677$ & $0.3379$ & $0.0435$ \\ \hline
        Star pyramid & Highly non-smooth & $72.5464$ & $2.4998$ & $0.1139$ \\ \hline
    \end{tabular}
\end{table}

\subsubsection{Multiply-connected domains}

The second plot from the left in Fig.~\ref{fig:eigen} presents the spectrum of the companion matrices for a multiply-connected toroidal surface with genus $g = 1$. Unlike the simply-connected case, the TD-PMCHWT spectrum for the torus contains only $N_p = 2 N_l - 4$ non-trivial eigenvalues around $1 + 0j$. The absence of four non-trivial eigenvalues can be attributed to the two global loops of the torus. Indeed, on a toroidal surface with genus $g$, the static outer and inner operators $\pm \tfrac{1}{2} Id - K_0$ each have a nullspace of dimension $g$, spanned by either poloidal or toroidal vector fields \cite{CAO+2009a}. Consequently, the TD-MFIOs $\KK$ and $\KK^\prime$ acting on these $2g$ global loops are of comparable magnitude to $\pm \tfrac{1}{2} Id$, displacing the associated $4g$ eigenvalues ($2g$ per unknown) away from $1+0j$ and into the interior of the unit circle.

This analysis reveals a deeper insight: on multiply-connected geometries, the TD-MFIOs strongly couple to global topological modes, modifying the spectral characteristics of the system.

\subsubsection{Highly non-smooth domains}

The effect of geometric properties on the spectral behavior of the MOT TD-PMCHWT system is further examined through a stability analysis on two highly non-smooth domains: the benchmark NASA almond and a star-based pyramid (see Fig.~\ref{fig:geometries}). These geometries are characterized by very sharp corners that significantly amplify quadrature errors in the discretization of the time-domain operators, thereby exacerbating instability (see Table~\ref{tab:shift_geometry} for further details). Moreover, the MFIOs on such non-smooth surfaces are no longer compact. Their contributions become non-negligible and can significantly shift the non-trivial eigenvalues of the TD-PMCHWT away from $1+0j$ to the inside of the unit circle. This behavior is confirmed by the two rightmost plots of Fig.~\ref{fig:eigen}, which illustrate a broad spread of the TD-PMCHWT eigenvalues from $1+0i$ to the interior of the circle.

\subsection{Effects of discretization parameters and material contrast}

\begin{figure}[!t]
    \centering
    \includegraphics[trim={0cm 0cm 0.3cm 0.7cm}, clip, width=0.99\linewidth]{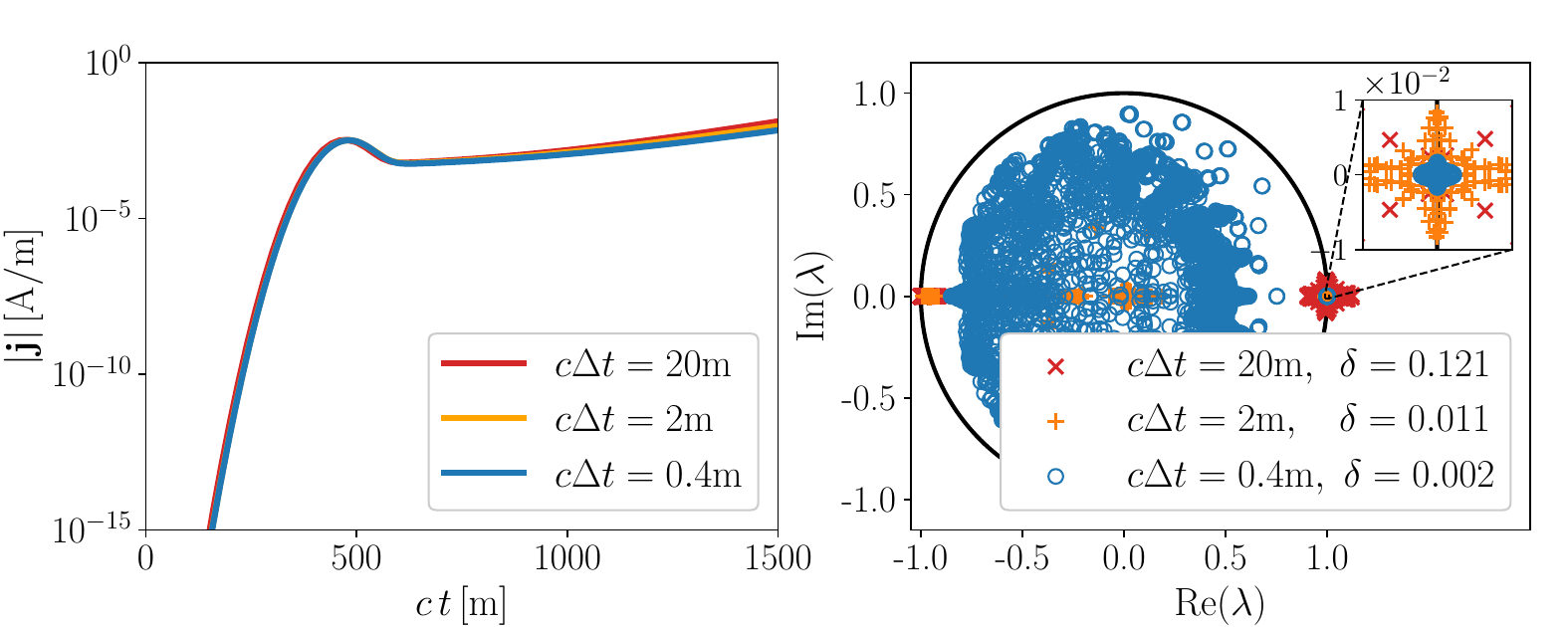}
    \caption{MOT solutions (\textit{left}) and spectrum of the companion matrices (\textit{right}) for the TD-PMCHWT on a sphere.
    The shift $\delta$ of eigenvalues near $1 + 0j$ increases linearly with $\Delta t$, but the growth of solution over time is independent.
    }
    \label{fig:eigen_dt}
\end{figure}

Fig.~\ref{fig:eigen_dt} illustrates the effect of the time step $\Delta t$ on the instability of the TD-PMCHWT solution. As $\Delta t$ increases, the shift $\delta$ of non-trivial eigenvalues near $1 + 0j$ increases linearly, while the MOT solutions exhibit similar exponential growth rates. This discrepancy does not contradict the statement in Section~\ref{sec:companion}. Indeed, MOT solution at late time $t_i = i \Delta t$ behaves as $\ub(t_i) \sim (1 + \delta)^i \approx \paren{1 + \delta/\Delta t}^{t_i}$. Thus, although $\delta$ linearly depends on $\Delta t$, the growth rate of MOT solution in continuous time is independent of it.

\begin{figure}[!t]
    \centering
    \includegraphics[trim={0.3cm 0cm 0.5cm 0.9cm}, clip, width=\linewidth]{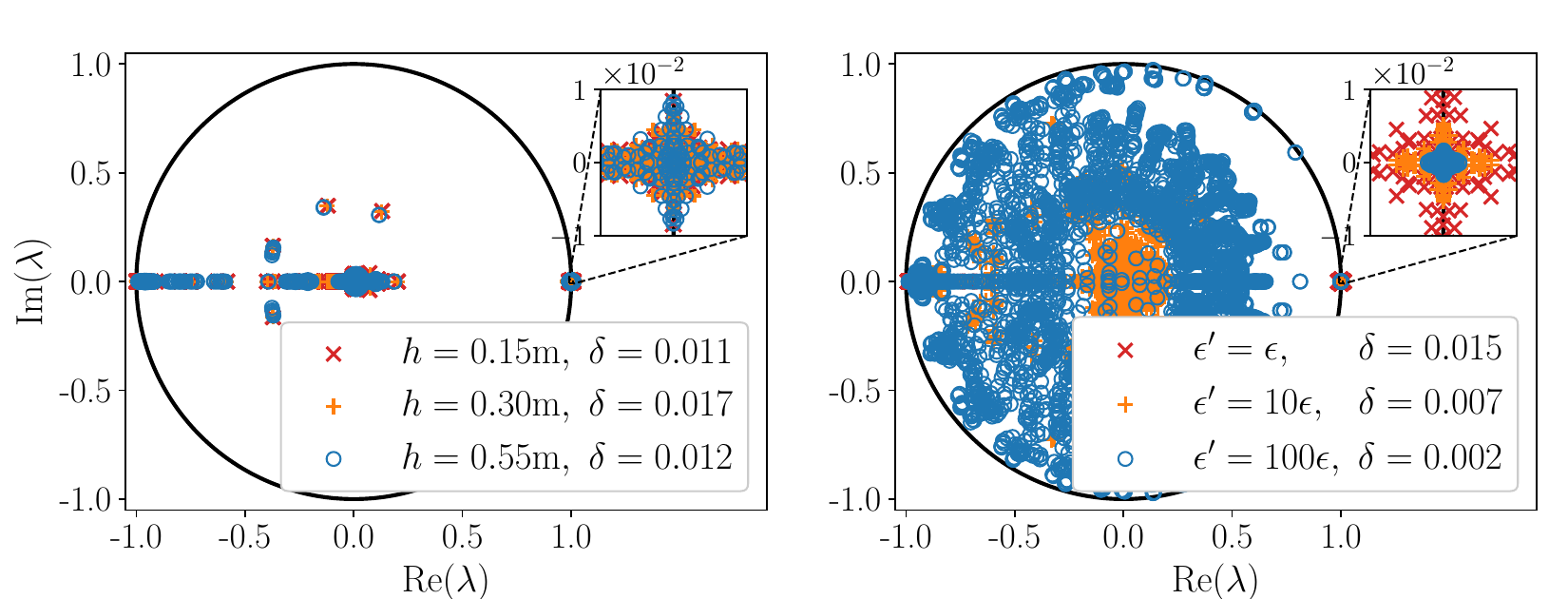}
    \caption{Spectrum of the companion matrices for the TD-PMCHWT on a sphere.
    (\textit{Left}) The surface is triangulated with different mesh sizes $h$.
    (\textit{Right}) 
    Different interior permittivities $\epsilon^\prime$ are used. The shift $\delta$ of eigenvalues around $1 + 0j$ is not significantly affected by $h$, while it decreases when increasing the material contrast.}
\label{fig:eigen_h}
\end{figure}

The influence of the mesh size $h$ and the material contrast are depicted in Fig.~\ref{fig:eigen_h}. Though decreasing the mesh size results in a quadratic increase of the number of spatial unknowns, the shift $\delta$ of non-trivial eigenvalues appears largely unaffected. In contrast, $\delta$ shows a decreasing dependence on the material contrast. This behavior can be explained by the fact that increasing the relative interior parameters $\epsilon^\prime$ and $\mu^\prime$ is effectively equivalent to decreasing the time step $\Delta t$ in the discretization of the interior operators $\TT^\prime$ and $\KK^\prime$.

\section{Conclusion}
\label{sec:conclusion}

The static solenoidal nullspace of the TD-EFIOs is the primary cause of the late-time instability of MOT solution to the TD-PMCHWT. The TD-MFIOs significantly influence the spectral characteristics of the MOT system, especially on multiply-connected surfaces or highly non-smooth domains. In all cases, numerical quadrature errors from integral evaluations have a substantial impact on the stability.

Building on these findings, future work will be aimed at resolving the late-time instability of the TD-PMCHWT solution by eliminating the nullspace of the TD-EFIOs. Some possible directions are outlined in \cite{BCA2015,BCA2015b,BCA2015d,LCA+2024}.

\ifCLASSOPTIONcaptionsoff
  \newpage
\fi

\IEEEtriggeratref{20}


\bibliographystyle{IEEEtran}
\bibliography{abrv_ref.bib}

%







\end{document}